\documentclass{article}

\usepackage{amssymb} 
\usepackage{amsmath} 
\usepackage{latexsym} 
\usepackage{theorem} 

\usepackage{epsfig}

\def\vect#1{\mbox{\boldmath $#1$}} 
\def\p{\prime}

\theorembodyfont{\itshape} 

\newtheorem{theorem}{Theorem}[section]
\newtheorem{lemma}[theorem]{Lemma}

\newtheorem{corollary}[theorem]{Corollary}
\newtheorem{proposition}[theorem]{Proposition}

\newtheorem{definition}[theorem]{Definition}

\newtheorem{conjecture}[theorem]{Conjecture}

\newenvironment{remark}{{\noindent \bf Remark:}}{\par\medskip}
\newenvironment{proof}{{\par\addvspace{0.1cm}\noindent \bf Proof. }}{\hfill$\Box$\par\medskip} 
 
\def\a{\alpha}
\def\w{\wedge}
\def\p{\prime}
\def\Si{\varSigma}
\def\Ga{\varGamma}
\def\R{\Re\mathfrak{e} \,}
\def\omegasub#1{\mbox{\large $\omega$}\mbox{\small${\mbox{\large ${}$}}_{#1}$}}
\def\spbmapright#1#2{\smash{%
 \mathop{\hbox to 0.8cm{\rightarrowfill}}
  \limits^{\displaystyle #1}_{\displaystyle #2}}}
\def\RR{\mathbb{R}}
\def\CC{\mathbb{C}}
\def\SS{\vect{S}}

\def\La{\Lambda}
\def\wb#1{\!\stackrel{\mbox{\boldmath$\underline{\mskip10mu}$}}{#1}}
\def\mb#1{\!\stackrel{\mbox{\boldmath$\underline{\mskip8mu}$}}{#1}}
\def\bb#1{\!\stackrel{\mbox{\boldmath$\underline{\mskip6mu}$}}{#1}}
\def\vect#1{\mbox{\boldmath $#1$}} 

\usepackage{color}
\def\pt#1{{\overset{\mbox{\large .}}{#1}}}

\title{\uppercase{Measure of a 2-component link}} %
\author{\textsc{Jun O'Hara}\footnote{ 
Partly supported by the Grant-in-Aid for Scientific Research (C),
Japan Society for the Promotion of Science. }
 }
\begin{document}

\maketitle


\begin{abstract}
A two-component link produces a torus as the product of the component knots in a two-point configuration space of a three-sphere. This space can be identified with a cotangent bundle and also with an indefinite Grassmannian. We show that the integration of the absolute value of the canonical symplectic form is equal to the area of the torus with respect to the pseudo-Riemannian structure, and that it attains the minimum only at the ``best'' Hopf links. 
\end{abstract}



\section{Introduction} 
Since energy of knots was introduced in \cite{OH1} about twenty years ago, aiming at producing an optimal knot for each knot type as an energy minimizer, a lot of related works have appeared, which form so-called {\sl geometric knot theory} (see, for example, \cite{CMR,CMRS,idealknots}). 
The present paper deals with the same type of topic. 
We introduce a functional on the space of $2$-component links such that the absolute minimum is attained only at ``best'' Hopf links, not at trivial links. 

Let $C_1\cup C_2$ be a $2$-component link in $\SS^3$. The value of our functional $A(C_1,C_2)$ can be interpreted in the following two ways. Observe that the link produces a torus $C_1\times C_2$ in $\SS^3\times\SS^3\setminus\Delta$, where $\Delta$ is the diagonal set. 

First, there is a natural identification between $\SS^3\times\SS^3\setminus\Delta$ and the total space of the cotangent bundle $T^\ast\SS^3$. The pull-back $\omega$ of the canonical symplectic form of $T^\ast\SS^3$ to $\SS^3\times\SS^3\setminus\Delta$ is the unique $2$-form (up to multiplication by a constant) which is invariant under the diagonal action of the M\"obius group. 
The $2$-form $\omega$ can also be considered as a natural symplectic form on the space of geodesics in a hyperbolic $4$-space $H^4$. 
As $\omega$ is exact, $\int_{C_1\times C_2}\omega$ vanishes, but $\int_{C_1\times C_2}|\omega|$ does not, which is $A(C_1,C_2)$. In this sense, it can be considered as an ``{\sl absolute symplectic measure}'' of the torus $C_1\cup C_2$ in $T^\ast\SS^3$. 

Second, from a M\"obius geometric viewpoint, $\SS^3\times\SS^3\setminus\Delta$ can be identified with the Grassmannian manifold $SO(4,1)/SO(1,1)\times SO(3)$ of oriented time-like $2$-dimensional vector subspaces in the $5$-dimensional Minkowski space $\RR^5_1$. 
By taking a pseudo-orthogonal complement of an oriented time-like $2$-dimensional vector subspace, we can identify this space with the Grassmannian manifold of oriented space-like $3$-dimensional vector subspaces in $\RR^5_1$. 

It has a natural pseudo-Riemannian structure which is compatible with the action of the Lorentz group, which induces the diagonal action of the M\"obius group to $\SS^3\times\SS^3\setminus\Delta$. Then $A(C_1,C_2)$ is equal to the measure (area) of the torus $C_1\cup C_2$ with respect to the pseudo-Riemannian metric. 

The key of the proof is that both the pull-back $\omega$ of the canonical symplectic form and the ``{\sl imaginary signed area element}'' with respect to the pseudo-Riemannian structure coincide with the real part of  the {\em infinitesimal cross ratio}, which \j{is a ``complex valued $2$-form'' on $C_1\times C_2$} used in the joint paper with Langevin \cite{La-OH1}. \j{Geometrically, it can be considered as} the cross ratio of $x, x+dx, y$ and $y+dy$, where these four points are considered as complex numbers by identifying a sphere through them with the Riemann sphere $\CC\cup\{\infty\}$. 

\j{Some remarks on the result on the energy of links \cite{AMN}, which is another characterization of the ``best'' Hopf link, 
will be given in Subsection \ref{AMN}.}

Throughout the paper, a link means a smooth (or at least of class $C^1$) $2$-component link. 

\medskip
{\bf Acknowledgment}. 
The author thanks deeply R\'emi Langevin, Masahiko Kanai and Luisa Paoluzzi for helpful suggestions. 
He also thanks the referee for a helpful suggestion concerning the symplectic form. 


\section{Two structures on {\boldmath $\SS^3\times \SS^3\setminus\Delta$}}\label{sec_S03}
We introduce two structures on $\SS^3\times \SS^3\setminus\Delta$, the symplectic structure and the pseudo-Riemannian structure, both compatible with M\"obius transformations. 
It is easy to see that both can be naturally generalized to $\SS^n\times \SS^n\setminus\Delta$ for any $n$. 

\subsection{Symplectic structure of $\SS^3\times \SS^3\setminus\Delta$}\label{subsect_sympl_str} 
\subsubsection{Via hyperbolic space}
As $\SS^3$ can be considered as the boundary of $4$-dimensional hyperbolic space $H^4$, $\SS^3\times \SS^3\setminus\Delta$ can be considered as the space of oriented geodesics in $H^4$, which is denoted by ${\mathcal{G}}$. 
The tangent space $T_\gamma {\mathcal{G}}$ along a geodesic $\gamma$ is the space of Jacobi fields along $\gamma$. 
Let $\nabla$ denote the Levi-Civita connection. 
Then, if we put 
$$\omega_g(\xi,\eta)=(\xi(t),\nabla_{\!\pt{\gamma}}\,\eta(t))-(\eta(t),\nabla_{\!\pt{\gamma}}\,\xi(t)) \hspace{0.5cm}(t\in\RR)$$
for $\xi,\eta\in T_\gamma {\mathcal{G}}$, where $(\>,\>)$ denotes the standard inner product on $T_{\gamma(t)} H^4$, then $\omega_g$ is an isometry-invariant symplectic form on ${\mathcal{G}}$ (see \cite[2C]{B}, \cite[3.1]{Kl}). 
Since an isometry of $H^4$ induces a M\"obius transformation of the boundary sphere $\SS^3$, $\omega_g$ defines a symplectic form on $\SS^3\times \SS^3\setminus\Delta$ which is invariant under the diagonal action of the M\"obius group. 

\subsubsection{Via cotangent bundle}
It is known that the space $\mathcal{G}$ of geodesics in $H^4$ is symplectomorphic to the cotangent bundle $T^\ast \SS^3$ \cite{F}. 
Let us give an identification between $\SS^3\times \SS^3\setminus\Delta$ and $T^\ast \SS^3$ explicitly. 

Assume $\SS^3$ is the unit sphere in $\mathbb{R}^{4}$. 
Let $\vect x$ be a point in $\SS^3$ and $p_{\mbox{\scriptsize \boldmath$x$}}\colon\SS^3\setminus\{\vect x\}\to(\textrm{Span}\langle\vect x\rangle)^\perp$ be a stereographic projection. By identifying $(\textrm{Span}\langle\vect x\rangle)^\perp$ with $T_{\mbox{\scriptsize \boldmath$x$}}\SS^3\cong T_{\mbox{\scriptsize \boldmath$x$}}^{\ast}\SS^3$, we obtain a bijection 
$$
\varphi_{\mbox{\scriptsize \boldmath$x$}}\colon\SS^3\setminus\{\vect x\}\ni\vect y\mapsto
\left(T_{\mbox{\scriptsize \boldmath$x$}}\SS^3\ni\vect v\mapsto 
p_{\mbox{\scriptsize \boldmath$x$}}(\vect y)\cdot \vect v
\in\mathbb{R}\right)\in T_{\mbox{\scriptsize \boldmath$x$}}^{\ast}\SS^3,
$$
where $\cdot$ denotes the standard inner product in $\mathbb R^{4}$. It induces a bijection 
\begin{equation}\label{phi_conf_sp_cot_bdl}
\varphi\colon
\SS^3\times \SS^3\setminus\Delta
\ni (\vect x,\vect y)
\,\mapsto\,
\left(\vect x, \varphi_{\mbox{\scriptsize \boldmath$x$}}(\vect y)\right)
\in T^{\ast}\SS^3.
\end{equation}

Let $\omegasub{S^3}$ be the canonical symplectic form of the cotangent bundle $T^{\ast}\SS^3$. Put $\omega=\varphi^\ast\omegasub{S^3}$. 
In \cite{La-OH1}, we showed that $\omega$ is invariant under the diagonal action of the M\"obius group. 
The converse is also true. Namely, if a $2$-form $\rho$ is invariant under the diagonal action of the M\"obius group, then $\rho=c\,\omega$ for some $c\in\RR$ (Proposition \ref{prop_kanai} in Appendix). 
Therefore, we can see that $\omega$ coincides with $\omega_g$ mentioned above up to a constant factor. 
%
\subsection{Pseudo-Riemannian structure of $\SS^3\times \SS^3\setminus\Delta$} 

The {\em Minkowski space} $\mathbb{R}^5_1$ is $\mathbb{R}^5$ with the indefinite inner product 
$$\langle \vect x, \vect y \rangle\!=\!-x_0y_0+x_1y_1+\cdots+x_4y_4.$$
The set of light-like vectors and the origin $\vect L=\left\{\vect v\in\mathbb{R}^{5}_1 \, ;\, \langle\vect v,\vect v\rangle=0\right\}$ is called the {\em light cone}. 
The $3$-sphere can be considered as the projectivization $\mathbb{P}\vect L$ of the light cone. It can also be identified isometrically with the intersection of the light cone and a  hyperplane given by $\{\vect x\,;\,\langle\vect x, \vect n\rangle=-1\}$, where $\vect n$ is a unit time-like vector. A $2$-dimensional vector subspace $\varPi$ of $\mathbb{R}^{5}_{1}$ is said to be time-like if $\langle\,,\,\rangle|_{\varPi}$ is non-degenerate and indefinite, namely, if $\varPi$ intersects the light cone transversely. 

A pair of points in $\SS^3$ can be considered as the intersection of $\SS^3$ and a $2$-dimensional time-like subspace of $\RR^5_1$. Therefore, if we also take the order of the points into account, $\SS^3\times \SS^3\setminus\Delta$ can be identified with the Grassmannian manifold $\widetilde{\textrm{\rm Gr}}_-(2;\mathbb{R}^{5}_1)$ 
of oriented $2$-dimensional time-like subspaces of $\mathbb{R}^{5}_1$, i.e., a homogeneous space $SO(4,1)/SO(3)\times SO(1,1).$ 

Let $\varPi$ be an oriented time-like $2$-dimensional plane spanned by an ordered basis $\{\vect u, \vect v\}$. 
Then $\varPi$ corresponds to a pure $2$-vector $\vect u\wedge\vect v\in\stackrel{2}{\mbox{$\bigwedge$}}\,\mathbb{R}^{5}_1$, which is determined by $\varPi$ up to a positive factor. 
As is stated on page 280 of \cite{hertrich}, $\vect u\w\vect v$ is time-like, i.e., $\langle \vect u\w\vect v, \vect u\w\vect v\rangle<0$, where the indefinite inner product on $\stackrel{2}{\mbox{$\bigwedge$}}\,\mathbb{R}^{5}_1$ is given by  
\[
\langle \vect u^1\wedge \vect u^2, \vect v^1\wedge  \vect v^2\rangle=\det\big(\langle \vect u^i, \vect v^j\rangle\big).
\]
On the other hand, it is known that a pure $2$-vector determines a $2$-plane. 
Thus the Grassmannian manifold $\widetilde{\textrm{\rm Gr}}_-(2;\mathbb{R}^{5}_1)$ of oriented $2$-dimensional time-like subspaces of $\mathbb{R}^{5}_1$ can be identified with the set of unit time-like pure $2$-vectors in $\stackrel{2}{\mbox{$\bigwedge$}}\,\mathbb{R}^{5}_1$, where the norm of $\stackrel{2}{\mbox{$\bigwedge$}}\,\mathbb{R}^{5}_1$ is given by $\|\vect v\|=\sqrt{|\langle\vect v, \vect v\rangle|}$. 
It is a $6$-dimensional pseudo-Riemannian manifold with index $3$. 
By taking a pseudo-orthogonal complement of an oriented time-like $2$-dimensional vector subspace, we can identify $\widetilde{\textrm{\rm Gr}}_-(2;\mathbb{R}^{5}_1)$ with the Grassmannian manifold $\widetilde{\textrm{\rm Gr}}_+(3;\mathbb{R}^{5}_1)$ of oriented space-like $3$-dimensional vector subspaces in $\RR^5_1$, which, in turn, can be identified with the set $\Theta(0,3)$ of unit space-like pure $3$-vectors in $\stackrel{3}{\mbox{$\bigwedge$}}\,\mathbb{R}^{5}_1$. 

Through the identifications mentioned above, the bijection from $\widetilde{\textrm{\rm Gr}}_-(2;\mathbb{R}^{5}_1)$ to $\widetilde{\textrm{\rm Gr}}_+(3;\mathbb{R}^{5}_1)$ is equal to the minus of the restriction of the Hodge $\star$ which is an isomorphism from $\stackrel{2}{\mbox{$\bigwedge$}}\,\mathbb{R}^{5}_1$ to $\stackrel{3}{\mbox{$\bigwedge$}}\,\mathbb{R}^{5}_1$ given by 
\[
\vect a\wedge \star\vect b=\langle\vect a, \vect b\rangle\, e_0\wedge e_1 \wedge \cdots \wedge e_4 \hspace{0.6cm}\big(\vect a,\vect b\in \stackrel{2}{\mbox{$\bigwedge$}}\,\mathbb{R}^{5}_1\,\big)
\]
(see \cite[p.288]{hertrich}). 

Let $\vect u$ and $\vect v$ be light-like vectors in $\mathbb R^5_1$. 
Put $\vect u\times \vect v=-\star(\vect u\wedge\vect v)\in \stackrel{3}{\mbox{$\bigwedge$}}\,\mathbb{R}^{5}_1$. 
Since the Hodge $\star$ satisfies $\langle \star \vect a, \star \vect b\rangle =-\langle \vect a, \vect b\rangle$, where $\vect a,\vect b\in \stackrel{2}{\mbox{$\bigwedge$}}\,\mathbb{R}^{5}_1$, we have 
\begin{equation}\label{f_<>_cross_prod}
\langle \vect u^1\times \vect u^2, \vect v^1\times  \vect v^2\rangle=-\det\big(\langle \vect u^i, \vect v^j\rangle\big).
\end{equation}

Thus we have a bijection 
\begin{equation}\label{def_psi}
\psi\colon\SS^3\times \SS^3\setminus\Delta\ni(\vect x, \vect y)\mapsto \frac{\vect x\times\vect y}{\,\|\vect x\times\vect y\|\,}\in\Theta(0,3).
\end{equation}
Since the indefinite inner product in \eqref{f_<>_cross_prod} is invariant under the action of the Lorentz group $O(4,1)$, the pseudo-Riemannian structure on $\SS^3\times \SS^3\setminus\Delta$ induced by $\psi$ is invariant under the diagonal action of the M\"obius group. 

\section{Measure of a 2-component link}\label{sec_main_thm}
All the pairs of points $\{(x,y)\,;\,x\in C_1, y\in C_2\}$ form a torus in $\SS^3\times \SS^3\setminus\Delta$. 
Let us call it the {\em product torus} of a $2$-component link $L=C_1\cup C_2$. 
%
\subsection{Area of the product torus of a link}
Let $\sigma$ be the composite of maps: 
\[\sigma\colon C_1\times C_2\stackrel{\iota}{\hookrightarrow} \SS^3\times \SS^3\setminus\Delta\,\spbmapright{\cong}{\psi}\,\Theta(0,3).\]
We identify $\sigma(C_1\times C_2)$ with $C_1\times C_2$ in what follows. 
The {\em area element} $dv$ of $C_1\times C_2$ associated with the pseudo-Riemannian structure of $\Theta(0,3)$ is given by 
$$
dv=\sqrt{\,\left|\det\left(\!\!\begin{array}{cc}
\langle \sigma_x, \sigma_x \rangle  &  \langle \sigma_x, \sigma_y \rangle \\
\langle \sigma_y, \sigma_x \rangle & \langle \sigma_y, \sigma_y \rangle 
\end{array}\!\!\right)\right|\,}\,dx\w dy\,,
$$
where $\sigma_x$ and $\sigma_y$ denote ${\partial \sigma}/{\partial x}(x,y)$ and ${\partial \sigma}/{\partial y}(x,y)$ in $T_{\sigma(x,y)}\Theta(0,3)$, respectively. 
\begin{definition} \rm 
Define the {\em measure} of a $2$-component link $L=C_1\cup C_2$ by the area of the product torus 
$$
A(C_1,C_2)=\int_{C_1\times C_2}dv
=\int_{C_1\times C_2}\sqrt{\,\left|\det\left(\!\!\begin{array}{cc}
\langle \sigma_x, \sigma_x \rangle  &  \langle \sigma_x, \sigma_y \rangle \\
\langle \sigma_y, \sigma_x \rangle & \langle \sigma_y, \sigma_y \rangle 
\end{array}\!\!\right)\right|\,}\,dx\w dy.
$$
\end{definition}
%
\subsection{Main Theorem}
%
%
\begin{theorem}\label{main_theorem}
\begin{enumerate}
\item The measure of a $2$-component link satisfies 
\begin{equation}\label{measure=symplectic}
A(C_1,C_2)=\int_{C_1\times C_2}|\iota^\ast\omega|,
\end{equation}
where $\iota$ is the inclusion from $C_1\times C_2$ into $\SS^3\times\SS^3\setminus\Delta$ and 
$\omega$ is the pull-back of the canonical symplectic form of $T^\ast\SS^3$ to $\SS^3\times\SS^3\setminus\Delta$.
\item The measure of a $2$-component link takes its minimum value $0$ if and only if $L$ is the image of the ``best'' Hopf link 
\begin{equation}\label{f_standard_Hopf_link}
\{(z,w)\in\CC^2;|z|=1, w=0\}\cup \{(z,w)\in\CC^2;z=0, |w|=1\}\subset\SS^3
\end{equation}
by a M\"obius transformation. 
\end{enumerate}
\end{theorem}

The equation \eqref{measure=symplectic} implies that the area of the product torus can also be called the ``{\sl absolute symplectic measure}'' of it. 

\smallskip
We prove the theorem in the next section. 
%
\subsection{Area element of a product torus in $\SS^3\times\SS^3\setminus\Delta$}\label{subsec_area_element}

\begin{lemma}\label{lem_null-vector_0-sphere}
Both $\sigma_x$ and $\sigma_x$ are null vectors, i.e., $\langle \sigma_x, \sigma_x \rangle=\langle \sigma_y, \sigma_y \rangle=0$. 
Therefore the area element $dv$ is given by $\sigma^{\ast}dv=|\langle \sigma_x, \sigma_y \rangle|\,dx\w dy.$
\end{lemma}

\begin{proof}
Suppose $\SS^3$ is embedded in $\RR^5_1$, and points in $C_1$ and $C_2$ are expressed by $\bar{x}(s)$ and $\bar{y}(t)$, respectively. Put $p(s,t)=\bar{x}(s)\times\bar{y}(t)$ and $\tilde\sigma(s,t)=\sigma(\bar x(s), \bar y(t))$. 
Then it is given by 
\[\tilde\sigma(s,t)=\displaystyle \frac{p(s,t)}{{\langle p(s,t), p(s,t)\rangle}^{1/2}}\,.\]
Since $\bar x$ and $\bar y$ are light-like vectors, the formula (\ref{f_<>_cross_prod}) implies  
%
$$
\langle {p}, {p}\rangle =\langle \bar{x}, \bar{y}\rangle^2, \>\>
\langle {p}, {p}_s\rangle = \langle \bar{x}, \bar{y}\rangle \langle \bar{x}_s, \bar{y}\rangle, \>\>
\langle {p}_s, {p}_s\rangle = \langle \bar{x}_s, \bar{y}\rangle^2.
$$
Therefore 
\[\langle \tilde\sigma_s, \tilde\sigma_s\rangle = \displaystyle \frac{\langle {p}, {p}\rangle \langle {p}_s, {p}_s\rangle-\langle {p}, {p}_s\rangle^2}{\langle {p}, {p}\rangle^2}=0.
\]
\end{proof}

We also put geometric explanation in Subsection \ref{subs_appendix_geom} in Appendix. 

\smallskip
Let us call  $\langle \sigma_x, \sigma_y \rangle\,dx\w dy$ the {\em imaginary signed area element} of a product torus $C_1\times C_2$.

\section{Proof of the Main Theorem}
%
\subsection{The infinitesimal cross ratio}
We assume that both components $C_1$ and $C_2$ are oriented. 
Suppose $x\in C_1$ and $y\in C_2$. Let $\Ga(x,x,y)$ be the circle which is tangent to $C_1$ at $x$ that passes through $y$, oriented by the tangent vectors to $C_1$ at $x$. Let $\theta$ $(0\le\theta\le\pi)$ be the angle between $\Ga(x,x,y)$ and the tangent vector to $C_2$ at $y$. 
We call it the {\em conformal angle} between $x$ and $y$ and denote it by $\theta_L(x,y)$. It was introduced by Doyle and Schramm. 

Let $\Omega_L$ be a complex valued $2$-form on $C_1\times C_2$ given by 
\begin{equation}\label{f_inf_cr}
\displaystyle \Omega_L(x,y)=e^{i\theta_L(x,y)}\frac{dx\w dy}{|x-y|^2}
\end{equation}
(see \cite{La-OH1}). 
As both the conformal angle $\theta_L$ and the $2$-form ${dxdy}/{|x-y|^2}$ are equivariant under the diagonal action of a M\"obius transformation $T$, so is $\Omega_{L}$, namely, $(T\times T)^{\ast}\Omega_{T(L)}=\Omega_{L}$ \cite{La-OH1}. 

Let us give a geometric interpretation of $\Omega_L$. 
Let $\Si_L(x,y)$ be a sphere that passes through four points $x, x+dx, y$ and $y+dy$, i.e., a sphere which is tangent to $C_1$ at $x$ and to $C_2$ at $y$. 
Let $p$ be a stereographic projection from $\Si_L(x,y)$ to $\mathbb{C}\cup\{\infty\}$ and $\tilde x$, $\tilde x+\widetilde{dx}$, $\tilde y$ and $\tilde y+\widetilde{dy}$ the images by $p$ of the four points $x, x+dx, y$ and $y+dy$, respectively. 
Then $\Omega_L(x,y)$ is equal to the cross ratio $(\tilde x+\widetilde{dx}, \tilde y ; \tilde x, \tilde y+\widetilde{dy})$: 
\begin{equation}\label{f_def_inf_cr}
\Omega_L(x,y)=\frac{\widetilde{dx}\widetilde{dy}}{(\tilde x-\tilde y)^2}
\sim\frac{(\tilde x+\widetilde{dx})-\tilde x}{(\tilde x+\widetilde{dx})
-(\tilde y+\widetilde{dy})}
:\frac{\tilde y-\tilde x}{\tilde y-(\tilde y+\widetilde{dy})}. 
\end{equation}
This is why we call $\Omega_L$ the {\em infinitesimal cross ratio}. 
We remark that the cross ratio does not depend on the stereographic projection $p$. 

\begin{remark}\label{remark31}\rm 
\j{The the form $dzdw/(z-w)^2$ on $\mathbb C\times \mathbb C\setminus\Delta$, which has been used in complex analysis, can also be obtained as the cross ratio of $w, w+dw, z$ and $z+dz$, as was mentioned by Rob Kusner, for example. 
In this sense, the infinitesimal cross ratio can be considered as generalization of $dzdw/(z-w)^2$ 
to a complex valued $2$-form on $C_1\times C_2$, or in general, $C\times C\setminus\Delta$, where $C$ is a union of space curves. 
In fact, when $C$ is a {\em plane} curve, the infinitesimal cross ratio can be obtained by restricting $dzdw/(z-w)^2$ to $C\times C\setminus\Delta$. In this case, it was used by H\'elein \cite{H} to show the isoperimetric inequality. }

\j{However, there is difficulty for {space} curves. 
First, $dzdw/(z-w)^2$ cannot be generalized to a $2$-form on the ambient space $\SS^3\times\SS^3\setminus\Delta$, so the restriction which works for the planar case does not work. To be precise, while the real part of $dzdw/(z-w)^2$ can be generalized to a $2$-form on $\SS^n\times\SS^n\setminus\Delta$ 
as we will see in the next subsection, the imaginary part 
cannot 
when $n\ge3$ as we will see in Proposition \ref{prop_kanai2}. }

\j{Secondly, even if we try to use the cross ratio to define the $2$-form, the cross ratio of four points in $\RR^n$ $(n\ge3)$ is not so well-behaved as in the planar case. This might be a reason why Ahlfors studied only the {absolute} cross ratio for the points in $\RR^n$ $(n\ge3)$ \cite{A}. 
When we want to define the cross ratio of (ordered) four points in $\RR^3$, we need the orientation of the sphere through the four points to avoid the ambiguity of complex conjugacy. 
There is a way to assign continuously the orientations to {all} the spheres given by the sets of ordered four points in $\RR^3$, i.e., there is a continuous map from ${\left(\RR^3\right)}^4\setminus\Delta$, where $\Delta$ is a big diagonal set, to the set of oriented $2$-speres in $\RR^3$, which can be identified with the de Sitter space in $5$-dimensional Minkowski space $\RR^5_1$. However, according to this method, the imaginary part of the cross ratio of any four points in $\RR^3$ is always non-negative (or, always non-positive according to the choice of a continuous map from ${\left(\RR^3\right)}^4\setminus\Delta$ to the de Sitter space). }
The reader is referred to \cite{OH2} for the details. 
\j{As a result, the imaginary part of the infinitesimal cross ratio may have singularity where it vanishes, just like that of the absolute value of a smooth function. Anyway, we do not use the imaginary part in this paper. }
\end{remark}

\subsection{The real part of the infinitesimal cross ratio}
In \cite{La-OH1} we showed that the pull-back of the canonical symplectic form of $T^\ast\SS^3$ to $C_1\times C_2$ coincides with the real part of the infinitesimal cross ratio up to a constant; 
\begin{equation}\label{real_part_of_inf_cr}
\iota^\ast\omega=\iota^{\ast}\varphi^{\ast}\omegasub{S^3}=-2\,\R\Omega_L
=-2\,\frac{\cos\theta_L(x,y)\,dx\w dy}{|x-y|^2}.
\end{equation}
\j{It seems that this fact in the case of $\SS^2$ is well known in symplectic geometry.}
\begin{lemma}\label{thm_Re_inf_cr=signed_area_element}
The imaginary signed area element of $C_1\times C_2$ with respect to the pseudo-Riemannian structure coincides with the real part of the infinitesimal cross ratio up to a constant; 
$$\langle \sigma_x, \sigma_y \rangle \, dx\wedge dy=2\,\R\Omega_{L}.$$
\end{lemma}

\begin{proof}
Suppose points in $C_1$ and $C_2$ are expressed as $x(s)$ and $y(t)$. Suppose $\SS^3$ is embedded in $\RR^5_1$ as the intersection of the light cone and a level hyperplane $\{x_0=1\}$. Let $\bar{x}$ and $\bar{y}$ be points in $\mathbb R^5_1$ corresponding to $x(s)$ and $y(t)$, i.e., $\bar{x}(s)=(1, x(s))$ and $\bar{y}(t)=(1, y(t))$. Put $\tilde\sigma(s,t)=\sigma(\bar x(s), \bar y(t))$ as before. 

The pull-back of the real part of the infinitesimal cross ratio is given by 
\begin{equation}\label{f1_re_infcr}
\left((x\times y)^{\ast}\R\Omega_{L}\right)(s,t)=\frac{\cos\theta_{L}(x(s),y(t))}{|x(s)-y(t)|^2}\,|x^{\p}(s)||y^{\p}(t)|\,ds\w dt.
\end{equation}
On the other hand, the pull-back of the imaginary signed area element is given by 
\begin{equation}\label{f2_signed_area_form}
\left((x\times y)^{\ast}\left(\langle \sigma_x, \sigma_y \rangle \, dx\w dy\right)\right)(s,t)=\langle\tilde\sigma_s, \tilde\sigma_t\rangle(s, t)\,ds\w dt.
\end{equation}
%

Fix any $(s_0, t_0)$. The M\"obius invariance of the both sides allows us to assume that $x(s_0)$ and $y(t_0)$ are antipodal. 
Then, at $(s_0, t_0)$, 
$$
\langle \bar x, \bar x\rangle=\langle \bar y, \bar y\rangle=0, \>\>
\langle \bar x, \bar y\rangle=-2. 
$$
Therefore, by the formula (\ref{f_<>_cross_prod}), at $(s_0, t_0)$ there holds 
$$
\langle p, p\rangle=4, \>\> \langle p, p_s\rangle=0, \>\> \langle p_s, p_t\rangle=-2x^{\p}(s_0)\cdot y^{\p}(t_0),
$$
which implies 
$$
\langle\tilde\sigma_s, \tilde\sigma_t\rangle=\frac{\langle p, p\rangle \langle p_s, p_t\rangle-\langle p, p_s\rangle \langle p, p_t\rangle}{\langle p, p\rangle^2}
=-\frac12\,x^{\p}(s_0)\cdot y^{\p}(t_0). 
$$
%
Since $x_0=x(s_0)$ and $y_0=y(t_0)$ are antipodal, we have (Figure \ref{pi-conf_angle})
\[\theta_L(x(s_0), y(t_0))=\pi-\angle x^{\p}(s_0)\cdot y^{\p}(t_0). \]
\begin{figure}[htbp]
\begin{center}
\includegraphics[width=.55\linewidth]{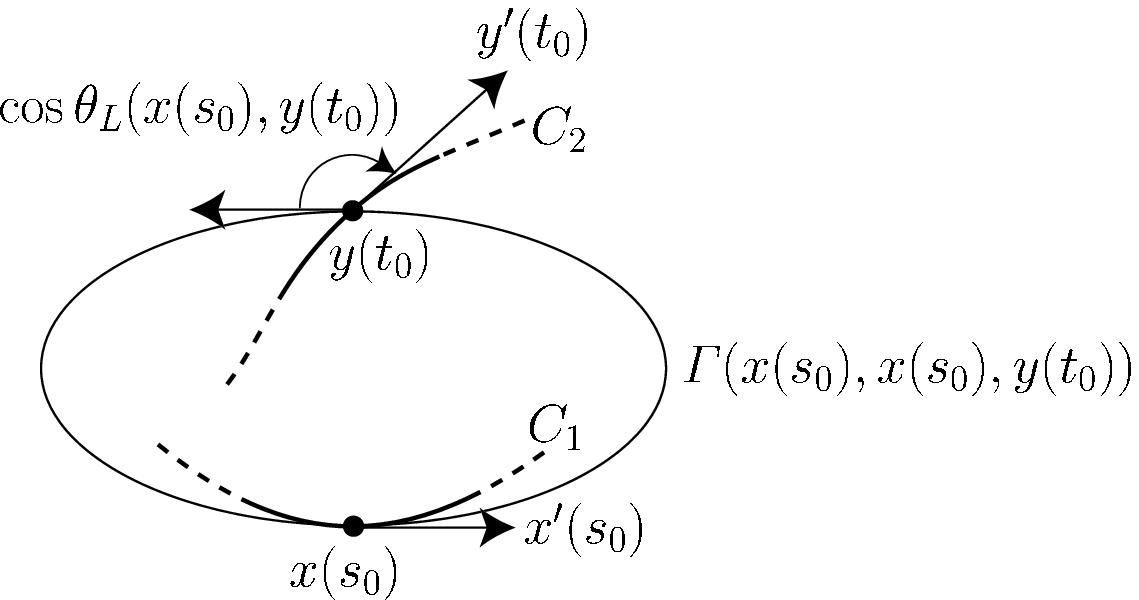}
\caption{}
\label{pi-conf_angle}
\end{center}
\end{figure}

\noindent
It follows that 
\[\langle\tilde\sigma_s, \tilde\sigma_t\rangle(s_0, t_0)=-\frac12\,x^{\p}(s_0)\cdot y^{\p}(t_0)
=2\frac{|x^{\p}(s_0)||y^{\p}(t_0)|}{|x(s_0)-y(t_0)|^2}\cos\theta_L(x(s_0), y(t_0)),\]
which implies that the right-hand sides of \eqref{f1_re_infcr} and \eqref{f2_signed_area_form} coincide. 
\end{proof}

We remark that an alternative geometric proof can be obtained if we use pseudo-orthonormal basis of $\SS^3\times \SS^3\setminus\Delta$ illustrated in Figure \ref{orthonormal_basis_S^3-2}. 
This is because 
$$
\langle \tilde\sigma_s+\tilde\sigma_t, \tilde\sigma_s+\tilde\sigma_t\rangle(s_0, t_0) =-x^{\p}(s_0)\cdot y^{\p}(t_0)
$$
implies 
$\langle\tilde\sigma_s, \tilde\sigma_t\rangle(s_0, t_0)=-(1/2)\,x^{\p}(s_0)\cdot y^{\p}(t_0).$

\begin{corollary}\label{cor_signed_area_form=sympl_form}
The imaginary signed area element of a product torus $C_1\times C_2$ with respect to the pseudo-Riemannian structure is equal to minus the pull-back of the canonical symplectic form:
$$
\langle \sigma_x, \sigma_y \rangle \, dx\wedge dy=-\iota^{\ast}\varphi^{\ast}\omegasub{S^3}.
$$
\end{corollary}
This completes the proof of Theorem \ref{main_theorem} (1). 

\medskip
We remark that a statement similar to that of the above corollary does not hold for a general surface in $\SS^3\times\SS^3\setminus\Delta$ as we will see in Subsection \ref{signed_area_form_not=sympl_form} in Appendix. 

\subsection{Proof of Theorem \ref{main_theorem} (2)}
As \setlength\arraycolsep{1pt}
\begin{equation}\label{f_area}
A(C_1,C_2)=2\int_{C_1\times C_2}\frac{|\cos\theta_L(x,y)|}{|x-y|^2}\,dx\,dy\,,
\end{equation}
it is equal to $0$ if and only if the conformal angle $\theta_L(x,y)$ is equal to ${\pi}/2$ for any $x\in C_1$ and $y\in C_2$. 

\begin{figure}[htbp]
\begin{center}
\begin{minipage}{.45\linewidth}
\begin{center}
\includegraphics[width=.75\linewidth]{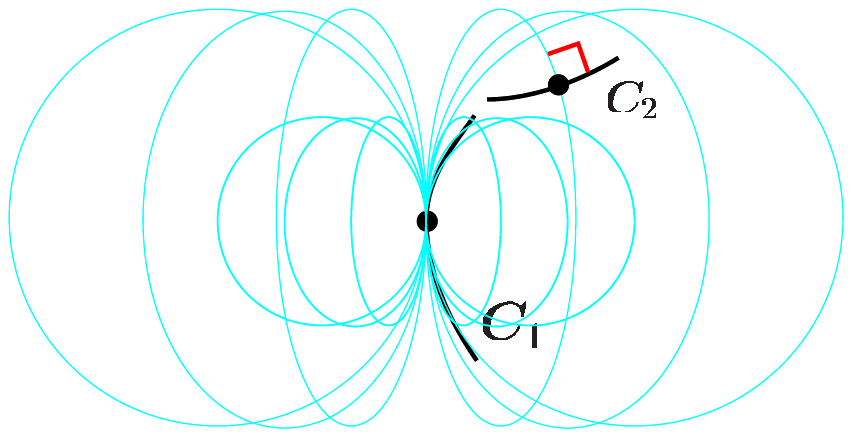}
\caption{The circles of $\mathcal{C}_{x}$}
\label{circles}
\end{center}
\end{minipage}
\hskip 0.4cm
\begin{minipage}{.45\linewidth}
\begin{center}
\includegraphics[width=.75\linewidth]{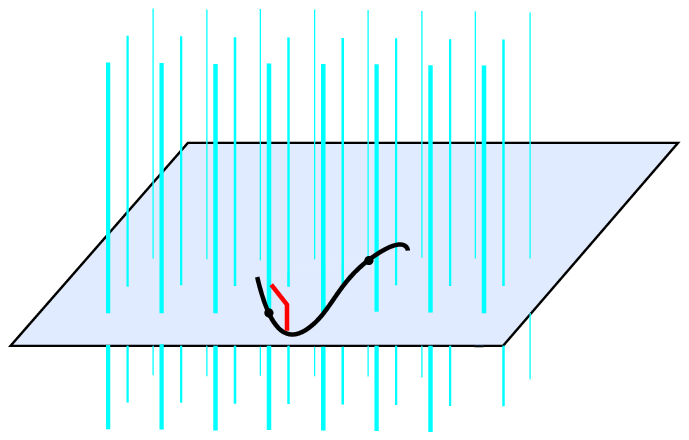}
\caption{The image by $\pi$}
\label{c1c2-2}
\end{center}
\end{minipage}
\end{center}
\end{figure}

\begin{figure}[htbp]
\begin{center}
\begin{minipage}{.3\linewidth}
\begin{center}
\includegraphics[width=.88\linewidth]{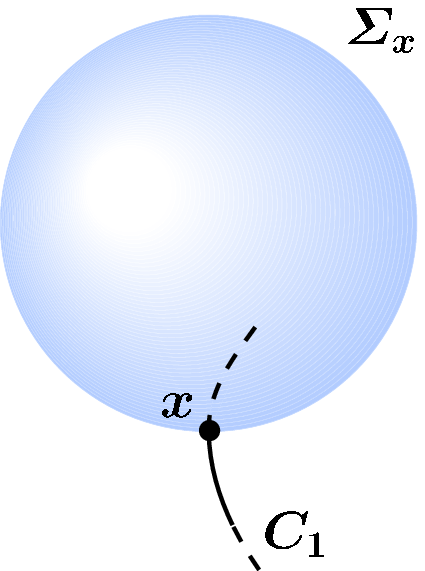}
\caption{}
\label{area_link_sphere}
\end{center}
\end{minipage}
\hskip 0.2cm
\begin{minipage}{.3\linewidth}
\begin{center}
\includegraphics[width=.88\linewidth]{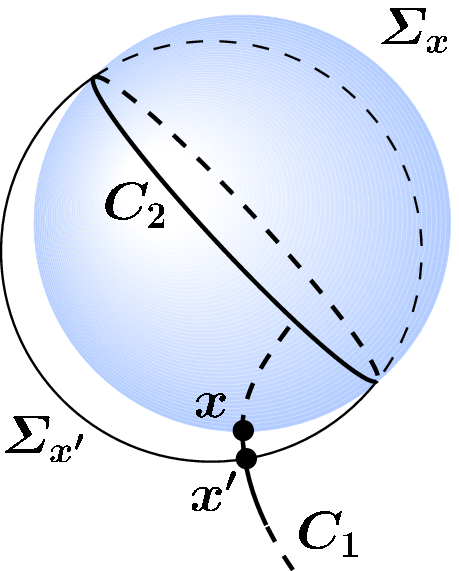}
\caption{}
\label{area_link_sphere2}
\end{center}
\end{minipage}
\hskip 0.2cm
\begin{minipage}{.34\linewidth}
\begin{center}
\includegraphics[width=.9\linewidth]{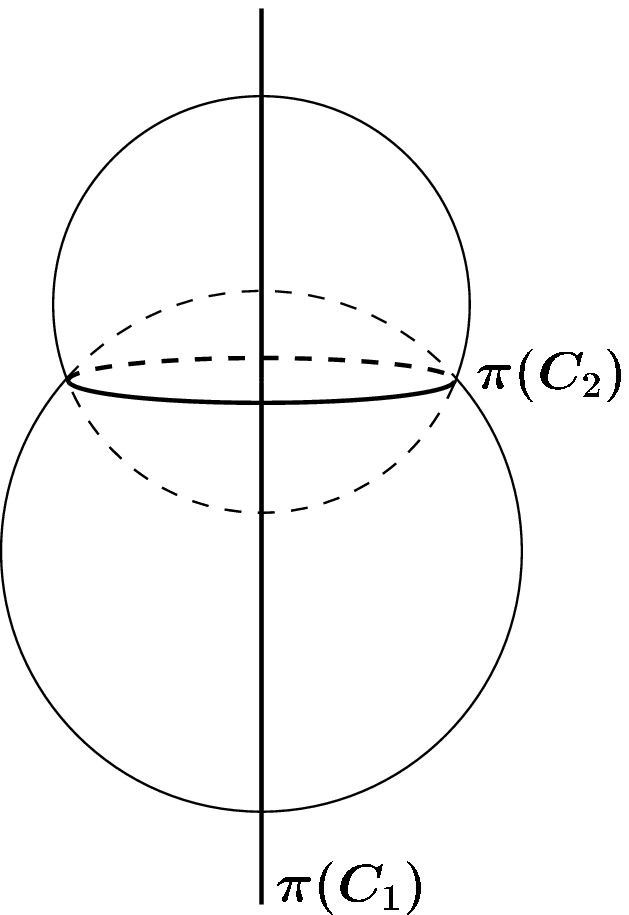}
\caption{}
\label{standard_Hopf_link}
\end{center}
\end{minipage}
\end{center}
\end{figure}

Suppose $A(C_1,C_2)=0$. 
Let $x$ be a point in $C_1$. 
Let $\mathcal{C}_{x}$ be the set of the circles which are tangent to $C_1$ at $x$. 
Then $C_2$ can intersect circles in $\mathcal{C}_{x}$ only at a right angle (Figure \ref{circles}). 
Consider a stereographic projection $\pi$ from $\SS^3\setminus\{x\}$ to $\RR^3$. 
It maps $\mathcal{C}_{x}$ to the set of parallel lines. 
Since $\pi(C_2)$ can intersect lines of $\pi(\mathcal{C}_{x})$ only at a right angle, $\pi(C_2)$ is contained in a $2$-plane which is orthogonal to the lines in $\pi(C_x)$ (Figure \ref{c1c2-2}). 
Therefore, $C_2$ is contained in a sphere $\Si_x$ which intersects $C_1$ at a right angle at $x$ (Figure \ref{area_link_sphere}). 

Let $x^{\p}$ be a point of $C_1$ close to $x$. 
As $C_1$ intersects $\Si_x$ orthogonally at $x$, we can take $x^{\p}$ outside $\Si_x$. 
Therefore, $\Si_x\ne\Si_{x^{\p}}$. 
Since $C_2$ is contained in the intersection $\Si_x\cap\Si_{x^{\p}}$, $C_2$ must be a circle (Figure \ref{area_link_sphere2}). 
The same argument shows that $C_1$ is also a circle. 

Consider the stereographic projection $\pi$ again. 
Since $C_1$ is a circle, $\pi(C_1)$ is a line. 
Then $\pi(C_2)$ is the intersection of two spheres which intersect the line $\pi(C_1)$ at a right angle (Figure \ref{standard_Hopf_link}). 
Therefore, $\pi(C_2)$ is symmetric in the line $\pi(C_1)$.  
It follows that $\pi(C_1)\cup \pi(C_2)$ is an image of the standard Hopf link (Figure \ref{standard_Hopf_link}). 

\medskip
This completes the proof of Theorem \ref{main_theorem} (2). 

\subsection{Corollary and Conjecture}
Let $[L]$ denote an isotopy class of a link $L$. Define 
\[A([L])=\inf_{C_1^{\p}\cup C_2^{\p}\in[L]}A(C_1^{\p},C_2^{\p}).\] 

\begin{corollary}
If $L$ is a separable link or a satellite link of a Hopf link, then $\mbox{\rm Area}\,([L])=0$. 
\end{corollary}

\begin{figure}[htbp]
\begin{center}
\includegraphics[width=.6\linewidth]{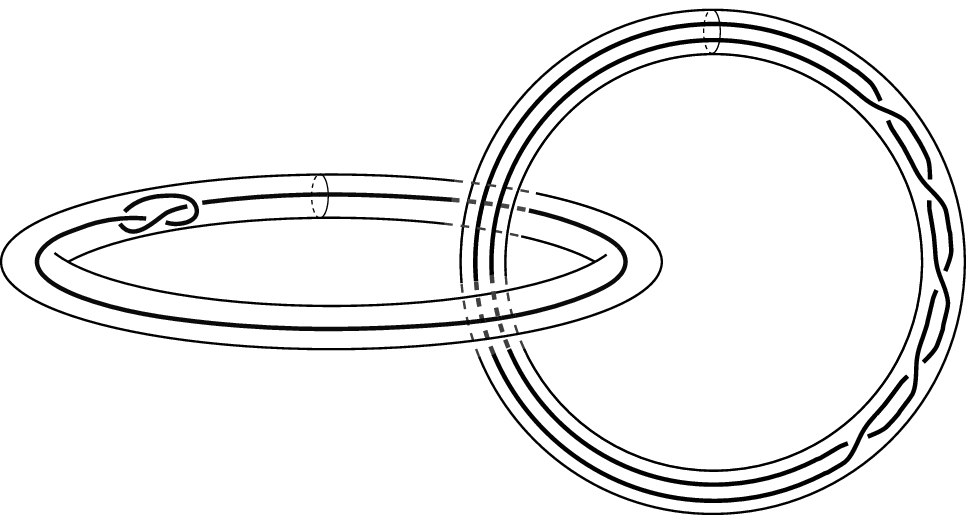}
\caption{A satellite link of a Hopf link}
\label{satellite_Hopf_link}
\end{center}
\end{figure}

\begin{proof}
Suppose $L=C_1\cup C_2$ is a separable link in $\RR^3$. 
We can make $|x-y|$ $(x\in C_1, y\in C_2)$ as big as we like. 
Now the conclusion follows from the formula (\ref{f_area}). 

\smallskip
Suppose $L=C_1\cup C_2$ is a satellite link of a Hopf link. 
Then, after an ambient isotopy, it can be contained in a very thin tubular neighbourhood of the standard Hopf link given by (\ref{f_standard_Hopf_link}). 
Furthermore, for any positive constants $\delta_1$ and $\delta_2$, the link can be placed so that, outside a small region of $C_1\times C_2$ whose measure is $\delta_1\textsl{Length}(C_1)\cdot\textsl{Length}(C_2)$, the conformal angle satisfies $|\theta_L-{\pi}/2|\le\delta_2$. 
Then the formula (\ref{f_area}) implies the assertion of the corollary since $|x-y|$ $(x\in C_1, y\in C_2)$ is bounded below. 
\end{proof}


\begin{conjecture} \rm 
We conjecture that $A([L])$ does not always vanish. 
For example, if $L=C_1\cup C_2$ is a hyperbolic link each component of which is a non-trivial knot, then there is no solid torus $H_1$ so that $C_1$ is contained in $H_1$ and $C_2$ in $\RR^2\setminus H_1$. 
We conjecture that $A([L])$ is positive for such a link type. 
\end{conjecture}
%
\section{Appendix}\label{appendix}
%
\subsection{Diagonal M\"obius invariance characterizes {$\vect{\omega}$}}\label{appendix1}
\begin{proposition}\label{prop_kanai}
Suppose $\rho$ is a $2$-form on $\SS^n\times \SS^n\setminus\Delta$ which is invariant under the diagonal action of M\"obius transformations. Then $\rho=c\,\omega$ for some constant $c$, where $\omega$ is the pull-back of the canonical symplectic form of $T^\ast \SS^n$ by the bijection from $\SS^n\times \SS^n\setminus\Delta$ to $T^\ast \SS^n$ given by \eqref{phi_conf_sp_cot_bdl}. 
\end{proposition}

This fact has been mentioned in \cite{Kanai} in a more general form (see $\S$ 3.2). We put the proof here since the author could not find it in the literature. 

\begin{proof}
%
We will show the equivalent statement for $\RR^n$. 
First note that the pull-back $\omegasub{{\RR^n}}$ of $\omega$ by a map $p^{-1}\times p^{-1}$ from $\RR^n\times \RR^n\setminus\Delta$ to $\SS^n\times \SS^n\setminus\Delta$, where $p$ is a stereographic projection from $\SS^n$ to $\RR^n$, is given by 
\begin{equation}\label{omegasub_R}
\begin{array}{rcl}
\omegasub{{\RR^n}}
%
&=&\displaystyle{2\left(\frac{\sum_{i=1}^{n} dX_i\w dY_i}{|\vect X-\vect Y|^2}
-2\frac{\big(\sum_{i=1}^{n} (X_i-Y_i)\,dX_i\big)\w 
\big(\sum_{j=1}^{n} (X_j-Y_j)\,dY_j\big)}{|\vect X-\vect Y|^4}\right)}
\end{array}
\end{equation}
 (see \cite{La-OH1}).

Let $\rho$ be a $2$-form on $\RR^n\times\RR^n\setminus\Delta$. Assume that $\rho$ is invariant under the diagonal action of M\"obius transformations. 
The invariance under the diagonal action of parallel translations implies that $\rho$ can be expressed as 
\begin{equation}\label{f_rho}
\rho=\sum_{i,j}f_{ij}(\vect x-\vect y)\,dx_i\w dy_j+\sum_{i<j}g_{ij}(\vect x-\vect y)\,dx_i\w dx_j+\sum_{i<j}h_{ij}(\vect x-\vect y)\,dy_i\w dy_j
\end{equation}
for some functions $f_{ij}, g_{ij}$ and $h_{ij}$. 

\smallskip
(i) Let us show $g_{ij}=h_{ij}=0$ $(i<j)$. Suppose $(i,j)=(1,2)$. 
The invariance under the diagonal action of rotation in $\mbox{\rm Span}\langle \vect e_1, \vect e_2\rangle$ 
shows 
\begin{equation}\label{219-g12-1}
\begin{array}{l}
g_{12}(v_1,\dots,v_n)
=\displaystyle g_{12}\Big(\sqrt{v_1^2+v_2^2}, \, 0, \,v_3,\dots,v_n\Big). \nonumber
\end{array}
\end{equation}
On the other hand, the invariance under the diagonal action of reflection in the hyperplane $\{v_2=0\}$ shows 
\begin{equation}\label{219-g12-2}
g_{12}(v_1,\dots,v_n)
=\displaystyle -g_{12}(v_1,-v_2,v_3\dots,v_n). \nonumber
\end{equation}
The above two equations 
imply $g_{12}=0$. 

\smallskip
(ii) Assume $n=2$. We show the statement using complex coordinates $w=x_1+ix_2$ and $z=y_1+iy_2$. The above argument shows that $\rho$ can be expressed as 
$$
\rho=F_1(w-z)\,dw\w dz+F_2(w-z)\,d\wb w\w d\bb z+F_3(w-z)\,d w\w d\bb z+F_4(w-z)\,d\wb w\w dz
$$
for some functions $F_i$. 
\begin{itemize}
%
\item[(a)] The invariance of $\rho$ under $(w,z)\mapsto(\beta w, \beta z)$ $(\beta\in\CC^\times)$ shows that 
$$
\left(\begin{array}{c}
F_1(\beta\zeta)\\
F_2(\beta\zeta)\\
F_3(\beta\zeta)\\
F_4(\beta\zeta)
\end{array}\right)
=\left(
\begin{array}{cc}
\begin{array}{cc}
\beta^{-2}&\\
&{\mb{\beta}}^{\,-2}
\end{array}
& \mbox{\Large\bf$0$} \\
\mbox{\Large\bf$0$}&
\begin{array}{cc}
|\beta|^{-2}&\\
&|\beta|^{-2}
\end{array}
\end{array}
\right)
\left(\begin{array}{c}
F_1(\zeta)\\
F_2(\zeta)\\
F_3(\zeta)\\
F_4(\zeta)
\end{array}\right).
$$
\item[(b)] The invariance of $\rho$ under $(w,z)\mapsto(\wb{w},\bb{z})$ shows that \setlength\arraycolsep{3pt}
\begin{equation}\label{F_bar}
\left(\begin{array}{c}
F_1\big(\mb{\zeta}\big)\\
F_2\big(\mb{\zeta}\big)\\
F_3\big(\mb{\zeta}\big)\\
F_4\big(\mb{\zeta}\big)
\end{array}\right)
=\left(
\begin{array}{cc}
\begin{array}{cc}
0&1\phantom{\mb{\zeta}}\hspace{-0.15cm}\\
1&0\phantom{\mb{\zeta}}\hspace{-0.15cm}
\end{array}
& \mbox{\Large\bf$0$} \\
\mbox{\Large\bf$0$}&
\begin{array}{cc}
0&1\phantom{\mb{\zeta}}\hspace{-0.15cm}\\
1&0\phantom{\mb{\zeta}}\hspace{-0.15cm}
\end{array}
\end{array}
\right)
\left(\begin{array}{c}
F_1(\zeta)\phantom{\mb{\zeta}}\hspace{-0.15cm}\\
F_2(\zeta)\phantom{\mb{\zeta}}\hspace{-0.15cm}\\
F_3(\zeta)\phantom{\mb{\zeta}}\hspace{-0.15cm}\\
F_4(\zeta)\phantom{\mb{\zeta}}\hspace{-0.15cm}
\end{array}\right).
\end{equation}
\setlength\arraycolsep{1pt}
\end{itemize}

Note that 
$$
\begin{array}{c}
\displaystyle |\zeta|^2F_3(\zeta)\stackrel{(\mbox{\scriptsize \j{a}})}=F_3(1)\stackrel{(\mbox{\scriptsize \j{b}})}=F_4(1)\stackrel{(\mbox{\scriptsize \j{a}})}=|\zeta|^2F_4(\zeta),\\[1mm]
\zeta^2F_1(\zeta)\stackrel{(\mbox{\scriptsize \j{a}})}=F_1(1)\stackrel{(\mbox{\scriptsize \j{b}})}=F_2(1)\stackrel{(\mbox{\scriptsize \j{a}})}={\bb{\zeta}}^{\,2}F_2(\zeta).
\end{array}
$$
Therefore, putting $a=F_1(1)$ and $b=F_3(1)$ we have 
$$
\rho=a\left(\frac{dw\w dz}{(w-z)^2}+\frac{d\wb w\w d\bb z}{(\wb w-\bb z)^2}\right)
+b\frac{d w\w d\bb z+d\wb w\w dz}{|w-z|^2}.
$$

\begin{itemize}
\item[(c)] Finally, the invariance of $\rho$ under $(w,z)\mapsto(1/w,1/z)$ implies $b=0$. 
\end{itemize}

Since 
$$
\frac{dw\w dz}{(w-z)^2}+\frac{d\wb w\w d\bb z}{(\wb w-\bb z)^2}=2\R \left(\frac{dw\w dz}{(w-z)^2}\right)
=-\omegasub{{\RR^2}},
$$
it completes the proof when $n=2$. 

\smallskip
(iii) Assume $n\ge3$. Put $c=-(1/2)f_{11}(\vect e_1)$, where $f_{11}$ is the function used in \eqref{f_rho} and $e_1$ is the first unit vector of $\RR^n$. 
If we apply the previous argument in (ii) to the $2$-planes $\textrm{Span}\langle \vect e_1, \vect e_j\rangle$ $(j\ne1)$ first and then to $\textrm{Span}\langle \vect e_i, \vect e_j\rangle$, we see that $\rho$ and $c\,\omegasub{\RR^n}$ have the same coefficients of $dx_i\w dy_j$ for all $i,j$, which completes the proof. 
%
\end{proof}

\j{We give another proposition which was announced in Remark \ref{remark31}. 
}
\j{
\begin{proposition}\label{prop_kanai2}
Suppose $\rho$ is a $2$-form on $\SS^n\times \SS^n\setminus\Delta$ that satisfies 
$(T\times T)^\ast \rho=(\textrm{\rm sgn}\, T) \rho,$ where $\textrm{\rm sgn}\, T$ is the signature of the Jacobian of $T$. 
Then $\rho=0$ if $n\ge3$, and $\rho=c\, \Im\mathfrak{m} \left(dz\wedge dw/(w-z)^2\right)$ for some constant $c$ if $n=2$ under the identification $\SS^2\cong \mathbb C\cup\{\infty\}$. 
\end{proposition}
}

\begin{proof}
\j{The proof goes somehow parallel to that of Proposition \ref{prop_kanai}. 
Suppose $\rho$ is expressed as \eqref{f_rho}. }

\smallskip
\j{(i) We can show $g_{ij}\equiv0$ and $h_{ij}\equiv0$ using the invariance up to sign under rotations and inversions in spheres. }

\smallskip
\j{(ii) Assume $n=2$. Then, the right-hand side of the formula \eqref{F_bar} is replaced by the minus of it, which implies that $\rho$ is a multiple of 
\[\frac{dw\w dz}{(w-z)^2}-\frac{d\wb w\w d\bb z}{(\wb w-\bb z)^2}=2\Im\frak{m} \left(\frac{dw\w dz}{(w-z)^2}\right).\]}

\smallskip
\j{(iii) Assume $n\ge3$. Let $(\vect x,\vect y)\in\RR^n\times\RR^n\setminus\Delta$, $\vect u\in T_x\mathbb R^n$, and $\vect v\in T_y\RR^n$. 
There is an orientation preserving M\"obius transformation $T$ such that $T(\vect x)=\vect 0, T(\vect y)=(1,0,\dots,0)$, and $T_\ast(\vect u)=\vect e_1$. Put $\widetilde{\vect v}=2\left(T_\ast(\vect v),\vect e_1\right)\vect e_1-T_\ast(\vect v)$. 
Then, as there is a rotation around $\vect e_1$-axis preserving $T(\vect x), T(\vect y)$ and $T_\ast(\vect u)$ that sends $\widetilde{\vect v}$ to $T_\ast(\vect v)$, we have $\rho(T_\ast(\vect u), T_\ast(\vect v))=\rho(T_\ast(\vect u), \widetilde{\vect v})$. 
On the other hand, as there is a reflection preserving $T(\vect x), T(\vect y)$ and $T_\ast(\vect u)$ that sends $\widetilde{\vect v}$ to $T_\ast(\vect v)$, we have $\rho(T_\ast(\vect u), T_\ast(\vect v))=-\rho(T_\ast(\vect u), \widetilde{\vect v})$. Hence $\rho\equiv0$. 
}
\end{proof}

\j{It follows that the imaginary part of $dz\wedge dw/(w-z)^2$ cannot be generalized to $\SS^n\times \SS^n\setminus\Delta$ when $n\ge3$. In fact, it can naturally be generalized to a K\"ahler form on $SO(n+1,1)/SO(2)\times SO(n-1,1)$, which is the space of oriented codimension $2$ spheres in $\SS^n$. }

\subsection{Pseudo-orthogonal basis of $\vect{\SS^3\times\SS^3\setminus\Delta}$}\label{subs_appendix_geom}
\begin{figure}[htbp]
\begin{center}
\includegraphics[width=.8\linewidth]{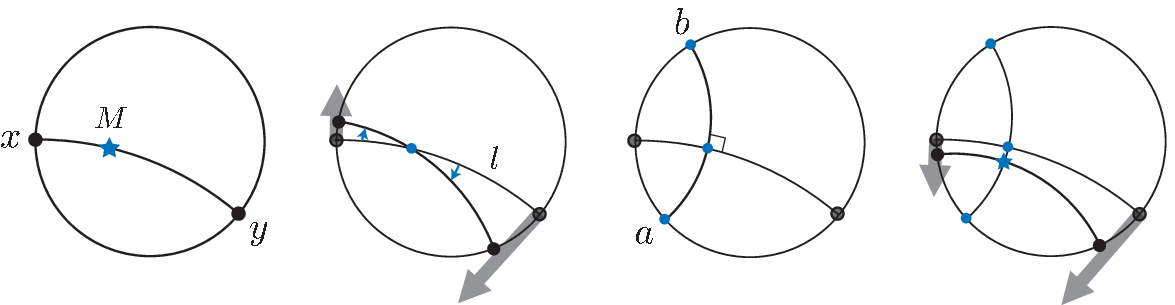}
\caption{Pseudo-orthogonal basis of a tangent space of $\SS^1\times\SS^1\setminus\Delta$. A picture in Poincar\'e disc model. }
\label{pencil_S^0_in_S^1}
\end{center}
\end{figure}

Let us start with a baby case $\SS^1\times\SS^1\setminus\Delta$. It can be identified with the set of oriented time-like planes in the $3$-dimensional Minkowski space $\RR^3_1$. By taking a positive unit normal vector to each of these planes, $\SS^1\times\SS^1\setminus\Delta$ can be identified with the $2$-dimensional de Sitter space $\La=\{\vect{x}\in\RR^3_1\,;\,\langle\vect x, \vect x\rangle=1\}$. 
Let $\Si=\{x,y\}$ be a pair of points in $\SS^1\cong\partial \mathbb{H}^2$. 
Let $l$ denote the geodesic in $\mathbb{H}^2$ which joins $x$ and $y$. 
Take a point $M$ on $l$ (Figure \ref{pencil_S^0_in_S^1}), then it determines two pencils as follows. 

Let $\vect a$ and $\vect b$ be the ``end points" of the geodesic in $\mathbb{H}^2$ which is orthogonal to $l$ at point $M$ (the third of Figure \ref{pencil_S^0_in_S^1}). 
Let $\mathcal{P}_+$ be a pencil obtained by rotating the geodesic $l$ around $M$ and $\mathcal{P}_-$ the Poncelet pencil with limit points $\vect a$ and $\vect b$. 
Then $\mathcal{P}_+$ and $\mathcal{P}_-$ can be considered as geodesics in $\La$, namely, the intersections with $\La$ and space-like and time-like $2$-planes $\Pi_\pm$. 
A pair of the unit tangent vectors to $\mathcal{P}_+$ and $\mathcal{P}_-$ at $\sigma$ can serve as a pseudo-orthonormal basis of $T_{\sigma}\La$, where $\sigma$ is a point in $\La$ that corresponds to $\Si$. 
These vectors can be obtained in $\Pi_\pm$ by rotation and Lorentz boost (hyperbolic rotation) of $\sigma$. 
The corresponding vectors in $\SS^1\times\SS^1\setminus\Delta$ are illustrated as the second and the last of Figure \ref{pencil_S^0_in_S^1}

\begin{figure}[htbp]
\begin{center}
\includegraphics[width=.5\linewidth]{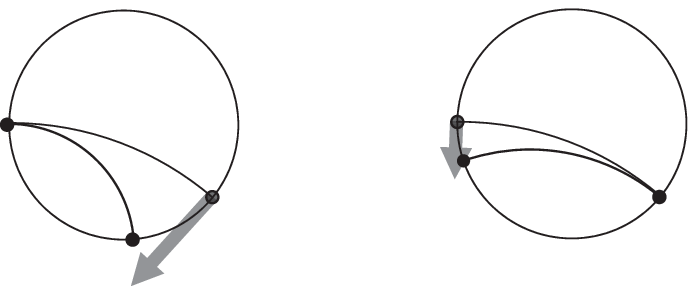}
\caption{Light-like basis of a tangent space of $\SS^1\times\SS^1\setminus\Delta$. A picture in Poincar\'e disc model. }
\label{pencil_S^0_in_S^1-l}
\end{center}
\end{figure}

Suppose $\{\vect u, \vect v\}$ is a pseudo-orthonormal basis of $T_{\sigma}\La$. Then we have another basis, $\left\{{(\vect u+\vect v)}/{\sqrt 2}, {(\vect u-\vect v)}/{\sqrt 2}\right\}$ consisting of two light-like vectors (Figure \ref{pencil_S^0_in_S^1-l}). 
This illustrates why $\sigma_x$ and $\sigma_y$ in Subsection \ref{subsec_area_element} are null vectors. 

\begin{figure}[htbp]
\begin{center}
\includegraphics[width=.66\linewidth]{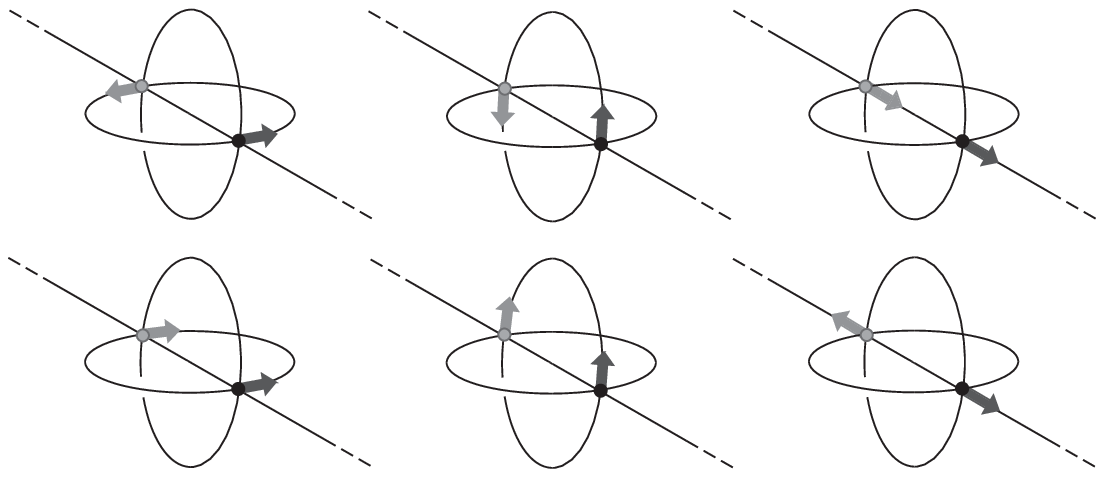}
\caption{Pseudo-orthogonal basis of a tangent space of $\SS^3\times\SS^3\setminus\Delta$. Space-like vectors above and time-like vectors below. A picture in $\mathbb{R}^3$ obtained through a stereographic projection}
\label{orthonormal_basis_S^3-2}
\end{center}
\end{figure}

The pseudo-orthonormal basis of ${\SS^3\times\SS^3\setminus\Delta}$ can be given by that of $\SS^1\times\SS^1\setminus\Delta$. In fact, we can consider three mutually orthogonal circles through a given pair of points, and take a pseudo-orthonormal basis in each circle as illustrated in Figure \ref{orthonormal_basis_S^3-2}. 
%

\subsection{The imaginary signed area element and the symplectic form}\label{signed_area_form_not=sympl_form}
Corollary \ref{cor_signed_area_form=sympl_form} does not necessarily hold for a surface in $\SS^3\times\SS^3\setminus\Delta$ which is not the product of two curves in $\SS^3$. 
Let us show it in $\RR^3\times\RR^3\setminus\Delta$, fixing a stereographic projection $p$ from $\SS^3$ to $\RR^3\cup\{\infty\}$. 

Suppose a pair of points in $\RR^3$ are expressed by $X(s,t)$ and $Y(s,t)$. 
Let $M$ be a surface $\{(X(s,t), Y(s,t))\}_{(s,t)\in D}$ in $\RR^3\times\RR^3\setminus\Delta$, where $D$ is a domain in $\mathbb R^2$. 
Put $X_s={\partial X}/{\partial s}, X_t={\partial X}/{\partial t}$, and 
$$
\widetilde X_s=2\left(X_s,\, \frac{X-Y}{|X-Y|}\right)\frac{X-Y}{|X-Y|}-X_s, \>
\widetilde X_t=2\left(X_t,\, \frac{X-Y}{|X-Y|}\right)\frac{X-Y}{|X-Y|}-X_t.
$$
Then $\widetilde X_s$ is the tangent vector at $Y$ to a circle which is tangent to $X_s$ at $X$ that passes through $Y$ with $\big|\widetilde X_s\big|=|X_s|$. The same interpretation also holds for $\widetilde X_t$. 

The pull-back of the canonical symplectic form of $T^\ast\RR^3\cong\RR^3\times\RR^3\setminus\Delta$ is given by 
$$
(X\times Y)^\ast\omegasub{\RR^n}=-2\left(\widetilde X_s\cdot Y_t-\widetilde X_t\cdot Y_s\right)\frac{ds\wedge dt}{|X-Y|^2}\,, 
$$
where $\omegasub{\RR^n}$ is given by \eqref{omegasub_R}. 
This can be verified by showing that the both sides coincide when $X$ and $Y$ are located on specific positions, say $X(s_0,t_0)=(1,0,0)$ and $Y(s_0,t_0)=(-1,0,0)$ because the both sides are equivariant under the diagonal action of M\"obius transformations. 

On the other hand, the ``{\sl signed area element}'' $\a_{\mbox{\tiny $M$}}$ of $M$ associated with the pseudo-Riemannian structure of $\Theta(0,3)$ can be given as follows. Let $\hat \sigma$ be the composite 
\[\hat \sigma:D\,
\smash{
 \mathop{\hbox to 1.2cm{\rightarrowfill}}
 \limits^{\displaystyle {X\times Y}}_{\displaystyle {}}}\,
M\hookrightarrow\RR^3\times \RR^3\setminus\Delta\,\spbmapright{p^{-1}\times p^{-1}}{}\,
\SS^3\times \SS^3\setminus\Delta\,\spbmapright{\cong}{\psi}\,\Theta(0,3).\]
Using the pseudo-orthonormal basis illustrated in Figure \ref{orthonormal_basis_S^3-2} and the M\"obius invariance, we have 
 \setlength\arraycolsep{1pt}
\[\begin{array}{rcl} 
(X\times Y)^\ast \a_{\mbox{\tiny $M$}}
&=&\displaystyle \sqrt{\,\det\left(\!\begin{array}{cc}
\langle \hat \sigma_s, \hat \sigma_s \rangle  \>&\>  \langle \hat \sigma_s, \hat \sigma_t \rangle \\
\langle \hat \sigma_t, \hat \sigma_s \rangle \>&\> \langle \hat \sigma_t, \hat \sigma_t \rangle 
\end{array}\!\right)}\,ds\w dt \\[4mm]
&=&\displaystyle 2\sqrt{\,\det\left(\!\begin{array}{cc}
2\widetilde X_s\cdot Y_s  \>\>&\>\>  \widetilde X_s\cdot Y_t+\widetilde X_t\cdot Y_s \\
\widetilde X_s\cdot Y_t+\widetilde X_t\cdot Y_s \>\>&\>\> 2\widetilde X_t\cdot Y_t
\end{array}\!\right)}\>\frac{ds\wedge dt}{|X-Y|^2}\,.
\end{array}\]

Therefore, the imaginary signed area element $\sqrt{-1}\,\a_M$ coincides with the pull-back of the canonical symplectic form $\omegasub{\RR^3}\big|_{C_1\times C_2}$ up to sign 
%
%
if and only if 
$
(\widetilde X_s\cdot Y_t)(\widetilde X_t\cdot Y_s)=(\widetilde X_s\cdot Y_s)(\widetilde X_t\cdot Y_t),
$
which holds if and only if $\widetilde X_s\times \widetilde X_t\perp Y_s\times Y_t.$ 
It does not always hold in general. 

\smallskip
We remark that this condition does not necessarily imply that the surface is a product of two curves. 
We also remark that the above condition is always satisfied for a surface in $\SS^1\times\SS^1\setminus\Delta$. 

\j{
\subsection{Remark on energy minimizing Hopf links}\label{AMN}
There is another variational characterization of the ``best'' Hopf link. }

\j{The {\em M\"obius cross energy} \cite{FHW} of a $2$-component link $C_1\cup C_2$, which is generalzation of the energy for knots defined by the author \cite{OH1}, is given by 
\[E(C_1,C_2)=\int_{C_1\times C_2}\frac{dxdy}{|x-y|^2}.\]
This energy is also invariant under M\"obius transformations. 
Recently, Agol, Marques and Neves proved  Freedman-He-Wang's conjecture, namely, they showed that if the linking number of $C_1$ and $C_2$ is equal to $\pm1$, then $E(C_1,C_2)\ge2\pi^2$, and that the equality holds if and only if $C_1\cup C_2$ is an image of the ``best'' Hopf link by a M\"obius transformation. 
This is a much more difficult problem, and was proved using min-max theory which has also been used in the proof of the Willmore conjecture \cite{MN}. }

\j{The formula \eqref{f_area} implies $E(C_1,C_2)\ge(1/2)A(C_1,C_2)$. To be more precise, the equality does not occur since the conformal angle between different components of a link cannot be identically zero. It might be interesting to point out that the infimum of $A(C_1,C_2)$ over all the $2$-component links is attained not at trivial links, but at the ``best'' Hopf link and the conformal image of it, whereas the infimum of $E(C_1,C_2)$ over all the $2$-component links is not attained, as $E(C_1,C_2)$ tends to $+0$ as the distance between $C_1$ and $C_2$ tends to $+\infty$. }


\bigskip
%

\par\noindent
Department of Mathematics, \\
Tokyo Metropolitan University, \\
1-1 Minami-Ohsawa, Hachiouji-Shi, \\
Tokyo 192-0397, Japan. \\

\noindent
{ohara@tmu.ac.jp}
%

\end{document}